\documentclass[12pt]{amsart}
\usepackage{amsfonts}
\usepackage{amssymb}
\usepackage{graphicx}

\def\ca{{\mathcal A}}

\def\Mp{M_{p\times p}}
\def\Np{N_{(p-1)\times (p-1)}}

\input{amssym.def}
\input{amssym.tex}
\newtheorem{theorem}{Theorem}[section]
\newcommand{\comment}[1]{}

\newtheorem{corollary}[theorem]{Corollary}
\newtheorem{remark}[theorem]{Remark}
\newtheorem{definition}[theorem]{Definition}

\title{Prime order automorphisms of Riemann surfaces}

\author{Jane Gilman}
\date{}
\thanks{This work was partially supported by a grants from the NSA and  the Rutgers
Research Council and by Yale University.}
\begin{document}

\begin{abstract}

Recently there has been renewed interest in the mapping-class group
of a compact surface of genus $g \ge 2$ and also in its finite order
elements. A finite order element of the mapping-class group will be
a conformal automorphisms on some Riemann surface of genus $g$. Here
we give the details of the proof that there is an {\sl adapted
basis} for any conformal automorphism of prime order on a surface of
genus $g$ and extend the original result to apply to fixed point
free automorphisms. An adapted basis is one that reflects the action
of the automorphism in the optimal manner described below. The proof
uses the Schreier-Reidemeister rewriting process. We find some new
consequences of the existence of an adapted basis. We also construct
an explicit example of such a basis and compute its intersection
matrix.
\end{abstract}
\maketitle
\section{Introduction} \label{section:intro}

 For a conformal automorphism of a
compact Riemann surface of the notion of an {\sl adapted} homology
bases was developed as part of a proof that the Riemann space (also
known as the Moduli Space) of  a punctured surface had the structure
of a quasi-projective variety \cite{Jdiss, Gth}. An adapted basis is
one which reflects the action of the  conformal automorphism in an
optimal way. Such an action would be reflected in the structure of
the period matrix of the surface in a useful manner. More recently
Rodriguez, Riera, Gonzalez and others  have used the notation of a
basis adapted to a group of automorphisms to obtain information
about abelian varieties and especially, the Prym variety. (See
\cite{RR}, \cite{GpRu}
 and the  references given there.)

 In this paper we survey earlier results about the matrix representation of a
 prime order automorphism with respect to an adapted basis and the corresponding
 intersection matrix for such a basis. The proof of the existence of an adapted basis uses the
Schreier-Reidemeister rewriting process. Here we give full details
of the application of the Schreier-Reidemeister rewriting process
used in \cite{G4} to construct the adapted basis. We have been told
that the application in \cite{G4} was too sketchy for some readers
to follow. We find some new consequences of the existence of an
adapted basis and extend the result to the case of a fixed point
free automorphism. We also construct an explicit example of such a
basis and compute its intersection matrix.

 The paper is organized as follows. In section
\ref{section:preliminaries} we fix notation and terminology and we
review the conjugacy invariants for an element of the mapping-class
group of prime order and basic facts about homology. Section
\ref{section:adapted} introduces the notion of an adapted homology
basis, section  \ref{section:exist} discusses the existence of such
bases and section \ref{section:intersection} the intersection
numbers of elements in an adapted basis. Section
\ref{section:matrixforms} fixes some matrix notation. In section
\ref{section:SCHR} the Schreier-Reidemeister rewriting process is
explained and the calculation is carried out in detail (section
\ref{section:rewrite}). The case for for a fixed point free
automorphism is carried out in section \ref{section:t=0} and some
corollaries are drawn in section \ref{section:additional}.

\tableofcontents

\section{Preliminaries} \label{section:preliminaries}
\subsection{Notation and Terminology}

 We let  $h$ be a conformal automorphism of
a compact Riemann surface $S$ of genus $g \ge 2$. Then $h$ will have
a finite number, $t$,  of fixed points. We let $S_0$ be the quotient
of $S$ under the action of the cyclic group generated by $h$ so that
$S_0 = S/\langle h \rangle$ and let $g_0$ be its genus. If $h$ is of
prime order $p$ with $p \ge 2$ , then the Riemann-Hurwitz relation
shows that
 $2g = 2pg_0 + (p-1)(t-2)$. If $p=2$, of course, this implies that $t$ will be
 even.

\subsection{Equivalent Languages} \label{section:tor}

 We emphasize that $h$ can be thought of in
a number or equivalent ways using different terminology. For a
compact Riemann surface of genus $g \ge 2$, homotopy classes of
homeomorphisms of surfaces are the same as isotopy classes.
Therefore, $h$ can be thought of as a representative of a homotopy
class or an isotopy class. Further, every isotopy class of finite
order contains an element of finite order so that $h$ can be thought
of as a homeomorphism of finite order. For every finite order
homeomorphism of a surface there is a Riemann surface on which its
action is conformal. A conformal homeomorphism of finite order up to
homotopy is finite. We use the language of conformal maps, but
observe that all of our results can be formulated using these other
classes of homeomorphisms.

We remind the reader that the Mapping-class group of a compact
surface of genus $g$ is also known as the Teichm{\"u}ller Modular
group or the Modular group, for short. We write $MCG(S)$ or
$MCG(S_g)$ for the mapping class group of the surface $S$ using the
$g$ when we need to emphasize that $S$ is a compact surface of genus
$g$. The Torelli Modular group or the Torelli group for short is
homeomorphisms of $S$ modulo those that induce the indentity on
homology and the homology of a surface is the abelianized homotopy.
There is surjective map $\pi$ from the mapping-class group onto
$Sp(2g,\mathbb{Z})$ that assigns to a homeomorphism the matrix of
its action on a canonical homology basis (see
\ref{section:homology}).

Since  $h$ can be thought of as a finite representative of a finite
order mapping-class, we will always treat it as finite.  For ease of
of exposition we use the language of a conformal maps and do not
distinguish between a homeomorphism that is of finite order or that
is of finite order up to homotopy or isotopy, a finite order
representative for the homotopy class, a conformal representative
for the class or the class itself. That is, between the topological
map, its homotopy class or a finite order representative or a
conformal representative.

For ease of exposition in what follows we first assume that  $t
> 0$. We treat the case $t=0$ separately in section
\ref{section:t=0}.

\subsection{Conjugacy Invariants for prime order mapping classes or
conformal automorphisms}

 Nielsen showed that the conjugacy class
of the image of $h$ in the mapping-class group is determined by a
set of  $t$ non-zero integers $\{n_1, ...., n_t\}$ with $0 < n_i <
p$ where $\Sigma_{i=1}^t n_i \equiv 0 \; (p)$. Here $\equiv \; (p)$
denotes equivalence modulo $p$.

 Let $m_j$ be the number of $n_i$
equal to $j$. Then we have $\Sigma_{i=1}^{p-1} i \cdot m_i \equiv 0
\; (p)$ (see \cite{G3} for details) and the conjugacy class is also
determined by the $(p-1)$-tuple, $(m_1,...,m_{p-1})$.

Topologically we can think of $h$ as a counterclockwise rotation by
an angle of  ${\frac{2\pi \cdot s_i}{p}}$ about the fixed point
$p_i, i = 1, ..., t$ of $h$. We call the $s_i$ the {\sl rotation
numbers}. The $n_i$ are the {\sl complimentary rotation numbers},
that is, $0 < s_i < p$ with $s_in_i \equiv 1 \; (p)$.

\subsection{Homology}\label{section:homology}
We recall the following facts about Riemann surfaces.

The homology group of a compact Riemann surface of genus $g$ is the
abelianized homotopy. Therefore, a homology basis for $S$ will
contain $2g$ homologously independent curves. Every surface has a
{\sl canonical homology basis}, a set of $2g$ simple closed curves,
$a_1,...,a_g; b_1,...,b_g$ with the property that for all $i$ and
$j$,  $a_i \times a_j = 0$,
 $b_i \times b_j = 0$ and $a_i \times
 b_j = \delta_{ij} = -b_j \times  a_i$ where $\delta_{ij}$ is the Kronecker delta.

\section{Adapted homology bases} \label{section:adapted}

 Roughly speaking a homology basis for $S$ is {\it adapted to $h$} if
it reflects the action of $h$ in a simple manner:  for each curve
$\gamma$ in the basis either all of the images of $\gamma$ under
powers of $h$ are also in the basis or the basis contains all but
one of the images of $\gamma$ under powers of $h$ and the omitted
curve is homologous to the negative of the sum of the images of
$\gamma$ under the other powers of $h$.

To be more precise

\begin{definition} \label{definition:adapt}  A homology basis
for $S$ is adapted to $h$ if for each $\gamma_0$ in the basis there
is a curve $\gamma$ with $\gamma_0 =h^k(\gamma)$ for  some integer
$k$ and either
\begin{enumerate}
\item  $\gamma, h(\gamma), ...h^{p-1}(\gamma)$ are all in the basis,
or
\item $\gamma, h(\gamma), ...h^{p-2}(\gamma)$ are all in the basis
and

$h^{p-1}(\gamma) \approx^h -(h(\gamma) +  h(\gamma)+ ...
+h^{p-2}(\gamma))$.

Here $\approx^h$ denotes is homologous to.

\end{enumerate}

\end{definition}

\section{Existence of adapted Homology Bases} \label{section:exist}
 It is known that
\begin{theorem} {\label{theorem:ad}} {\rm (Gilman 1977)} \cite{G4} There is a homology basis adapted to $h$. In
particular, if  $g\ge 2$, $t \ge 2$,  $g_0$ are as above, then
the adapted basis has $2p \times  g_0$ elements of type (1) above
and $(p-1)(t-2)$ elements of type (2).
\end{theorem}

and thus it follows that
\begin{corollary} \cite{G4} \label{corollary:corad}
Let  $M_{\ca}(h)$ denote the {\sl adapted matrix of $h$}, the
matrix of the action of $h$  with respect to an adapted basis. Then
$M_{\ca}(h)$ will be composed of diagonal blocks, $2g_0$ of which
are $p \times p$ permutation matrices with $1$'s along the super
diagonal and $1$ in the leftmost  entry of the last row and $t$ are
$(p-1) \times ( p-1)$ matrices with $1$'s along the super diagonal
and all entries in the last row $-1$.
\end{corollary}

A proof of theorem \ref{theorem:ad} is given in section
\ref{section:rewrite}.

\begin{remark} We adopt the following convention.
When we pass from homotopy to homology, we use the same notation for
the homology class of the curve as for the curve or its homotopy
class, but write $\approx^h$ instead of $=$. It will be clear from
the context which we mean.
\end{remark}

\section{Intersection Matrix for an Adapted Homology Basis}
\label{section:intersection}
%\section{The intersection matrix for an adapted homology basis}
So far  information about $M_{\ca}(h)$ seems to depend only on $t$
and not upon the $(p-1)$-tuple  $(m_1,...,m_{p-1})$ or equivalently,
upon the set of integers $\{n_1, ..., n_t\})$ which determines
 the conjugacy class of $h$ in the mapping-class group.  However,
while the $2pg_0$ curves can be extended to a canonical homology
basis for $h$, the rest of the basis can  not and its intersection
matrix, $I_{\ca}$ depends upon these integers.

In \cite{GP} the intersection matrix for the adapted basis was
computed.

The adapted basis consisted of the curves of type (1):
$$\{A_w, B_w,  w=1,..., g_0 \} \cup \{ h^j(A_w), h^j(B_w),
j=1,...,p-1 \}$$ and (some of) the curves of type (2):
$$X_{i,v_i}, h^j(X_{i,v_i}),\;  i=1,...,(p-1),\; j =1, ..., p-2,\;
v_i=1,...,u_i.$$

\begin{remark} For any one reading the original paper \cite{GP} note that the roles of $m$
and $n$ are interchanged here. To avoid confusion, we use $u$ and
$v$ in this section.
\end{remark}

%reversed there frHere we let $u=n$ and $v=m$ however, we note that
%in the original paper the roles of $m$ and $n$ are switched from the
%notation here. \marginpar{\bf HERE $n_i$ and $m_i$ have been
%switched!!!}

%The symbol $X_{i,m_i)$ represents a curves with

A lexicographical order is placed on $X_{i,v_i}$ so that $(r,v_r) <
(s, v_s)$ if and only if $r < s$ or $r=s$ and $v_r < v_s$. The $t-2$
curves $X_{s,v_s}$ with the largest subscript pairs are to be
included in the homology basis. Let $\hat{s}$ be the smallest
integer $s$ such that $u_s \ne 0$ and let $\hat{q}$ be chosen so
that $\hat{q} \hat{s} \equiv 1 (p)$. For any integer $v$ let $[v]$
denote the least non-negative residue of $\hat{q} v$ modulo $p$.
Thus the integer $[v]$ satisfies $0 \le [v] \le p-1$ and $\hat{s}
\times [v] \equiv v $ mod(p).

\comment{For an explanation of why we  normalize to $[v]$ rather
than using $v$ see section \ref{section:m[m]}.}

\begin{theorem} {\label{theorem:GPint}} {\rm (Gilman-Patterson,
1981)}
\cite{GP}  If $(u_1,...,u_{p-1})$ determines the conjugacy class of
$h$ in the mapping-class group, then the surface $S$ has a homology
basis consisting of:

\begin{enumerate}
\item  $h^j(A_w), h^j(B_w)$  where $1 \le w \le g_0,  0 \le j \le p-1$.

\item $h^k(X_{s,{v_s}})$   where $0 \le k \le p-2$ and for all pairs $(s,v_s)$
with $1 \le s \le p-1, 1 \le v_s \le u_s$ except that the two
smallest pairs are omitted.
%\end{itemize}

 The intersection numbers for the elements of the
adapted basis are given by
\begin{enumerate}
\item $h^j(A_w) \times h^j(B_w) = 1$

\item If $(r,v_r) < (s,v_s)$, then

$ {h^0(X_{r,{v_r}}) \times h^k(X_{s,{v_s}})} = \left\{
\begin{array}{ll}
        1 & {\mbox{ if }} \;\; [k] < [r] \le [k+s]    \\
       -1 & \mbox{ if } \;\; [k+s] < [r] \le [k]
                                    \end{array}
                                         \right. $

$ {h^0(X_{s,{v_s}}) \times h^k(X_{s,{v_s}})}
                  =  \left\{  \begin{array}{ll}
                        1 & {\mbox{ if } \;\; [k] \le [s] <  [k+s]} \\
\                       -1 & {\mbox{ if } \;\; [k+s] < [s]  < [k]}
                   \end{array}
                      \right.  $
\end{enumerate}

\item All other intersection numbers are $0$ except for those that
following from the above by applying the identities below to arbitrary
homology classes $C$ and $D$.
\subitem $ C \times D = - D \times C$

\subitem $ h^j(C) \times h^k(D)  = h^0(C) \times h^{k-j}(D)$, ($k-j$
reduced modulo $p$.)
\end{enumerate}
\end{theorem}

\begin{proof} For details we refer the reader to \cite{GP}.
Basically, the proof of this theorem comes from a careful
interpretation of the isomorphism between covering groups,
fundamental groups, and defining subgroups of coverings and their
relation to words corresponding to closed curves on the quotient
surface that lift to closed curves.
\end{proof}

\section{Matrix forms \label{section:matrixforms}}
%\marginpar{move this to later}

 We can write the results of theorems
\ref{theorem:ad} and \ref{theorem:GPint} and corollary
\ref{corollary:corad} in an explicit matrix form. To do so we fix
notation for some matrices. We will use the various explicit forms
in subsequent sections.

We let $M_{{\tilde{\mathcal{A}}}}$ denote the matrix of the action
of $h$ on an adapted basis and $I_{{\tilde{\mathcal{A}}}}$ be the
corresponding intersection matrix. Further, we let $M_{h_{CAN}}$ be
the matrix of the action of $h$ on a canonical homology basis. The
corresponding intersection matrix is denoted by $I_{h_{CAN}}$. If we
let $I_k$ denote the $k\times k$ identity matrix, then $I_{h_{CAN}}$
is (conjugate to) the $2g\times2g$ matrix $ = \left(
\begin{array}{cc}
0 & I_g\\
-I_g & 0 \\
\end{array} \right). $

However, we prefer to replace it by the following block matrix where
$q= {\frac{(p-1)(t-2)}{2}}$

$$I_{h_{CAN}} = \left(
\begin{array}{cccc}
0 & I_{pg_0} & 0 & 0\\
-I_{pg_0} & 0 & 0 & 0\\
0 & 0 & 0 & I_{q} \\
0 & 0 & -I_q & 0 \\
\end{array} \right). $$

 We denote the
$p\times p$ permutation matrix by

$$M_{p\times p}  =  \left(
\begin{array}{ccccccc}
0 & 1 & 0 & 0  &\ldots &0 & 0\\
0 & 0 & 1 & 0  &\ldots &0 & 0\\
0 & 0 & 0 & 1  &\ldots &0&  0\\
& & & & \ddots \\
%\vdots & \vdots&\vdots&\vdots&\ldots &\vdots&\vdots\\
0 & 0 & 0 & 0  &\ldots &1& 0\\0 & 0 & 0 & 0  &\ldots& 0 & 1\\
1 & 0 & 0 & 0  &\ldots &0 & 0\\
\end{array} \right);$$

the $(p-1)\times (p-1)$ non-permutation matrix of the theorem by

$$N_{(p-1)\times(p-1)}  = \left(
\begin{array}{rrrrrrr} %lllllll} %ccccccc}
0 & 1 & 0 & 0  &\ldots &0 & 0\\
0 & 0 & 1 & 0  &\ldots &0 & 0\\
0 & 0 & 0 & 1  &\ldots &0&  0\\
& & & & \ddots \\
%\vdots & \vdots&\vdots&\vdots& \ldots&\vdots&\vdots\\
0 & 0 & 0 & 0  &\ldots &1& 0\\0 & 0 & 0 & 0 & \ldots &0 & 1\\
-1 & -1 & -1 & -1&  \ldots& -1 & -1\\
\end{array} \right) $$

%$\Mp$, $\Mpt$

Thus we have the  $2g_0 \times p^2$ block matrix

$$M_{{\mathcal{A}}_{2g_0,  p\times p}} = \left(
\begin{array}{ccccccc}
\Mp & 0& 0  &\ldots &0 & 0\\
0 & \Mp & 0   &\ldots &0 & 0\\
& & & \ddots \\
%\vdots & \vdots&\vdots&\vdots&\ldots &\vdots&\vdots\\
0 & 0 &  0  &\ldots &\Mp& 0
\\0 & 0  & 0  &\ldots& 0 & \Mp
\\
\end{array} \right) $$

and the $(t-2)\times (p-1)^2$ block matrix

 $$N_{{\mathcal{A}}_{(t-2), (p-1)\times (p-1)}} = \left(
\begin{array}{ccccccc}
\Np & 0& 0   &\ldots &0 & 0\\
0 & \Np &  0  &\ldots &0 & 0\\
& & & \ddots& \\
%\vdots & \vdots&\vdots&\vdots&\ldots &\vdots&\vdots\\
0 & 0 & 0   &\ldots &\Np& 0
\\0 & 0 &  0  &\ldots& 0 & \Np\\
\end{array} \right) $$ so that the $2g \times 2g $ matrix $M_{\mathcal{A}}$ breaks into
blocks and can be written as
$$M_{\mathcal{A}} = \left(
\begin{array}{cc}
M_{{\mathcal{A}}_{2g_0, p\times p}} &  0 \\
%0 & M_{{\mathcal{A}}_{g_0, p\times p}} &  0 \\
 0&  N_{{\mathcal{A}}_{(t-2),  (p-1)\times (p-1)}}
\\
\end{array} \right) $$ where the blocks are of appropriate size.
The basis can be rearranged so that  $2g \times 2g $ matrix
$M_{\mathcal{\tilde{A}}}$ corresponding to the rearranged basis
breaks into blocks
$$M_{\mathcal{\tilde{A}}} = \left(
\begin{array}{ccc}
M_{{\mathcal{A}}_{g_0, p\times p}} &  0 & 0\\
0 & M_{{\mathcal{A}}_{g_0, p\times p}} &  0 \\
0 & 0&  N_{{\mathcal{A}}_{(t-2),  (p-1)\times (p-1)}}
\\
\end{array} \right) $$ Here
 the submatrix $$\left( \begin{array}{cc}
M_{{\mathcal{A}}_{g_0, p\times p}} &  0 \\
0 & M_{{\mathcal{A}}_{g_0, p\times p}}  \\
%0 & 0&  N_{{\mathcal{A}}_{(t-2),  (p-1)\times (p-1)}}
%\\
\end{array} \right) $$ is a symplectic matrix. We obtain the corollary

\begin{corollary} Let $S$ be a compact Riemann surface of genus $g$ and
assume that $S$ has a conformal automorphism $h$ of prime order $p
\ge 2 $. Assume that $h$ has $t$ fixed points where $t \ge 2$. Let
$S_0$ be the quotient surface $S_0 = S / \langle h \rangle $ where
$\langle h \rangle$ denotes the cyclic group generated by $h$ and
let $g_0$ be the genus of $S_0$ so that  $2g = 2pg_0 + t(p-1)$.
\vskip .1in
 There is a homology bases on which the action of $h$
is given by the $2g \times 2g$ matrix $M_{\mathcal{\tilde{A}}}$.
\vskip .1in
 The matrix $M_{\mathcal{\tilde{A}}}$ contains a $2g_0p
\times 2g_0p$ symplectic submatrix, but is not a symplectic matrix
except in the special case $t = 2$.
\end{corollary}

\begin{remark} We note that if $p=2$,
$\Mp$ reduces to $ \left(
\begin{array}{cc}
0& 1\\
1 &  0
\\
\end{array} \right) $ and $\Np$ to the $1\times1$ matrix $-1$.
\end{remark}

The point here is that while two automorphisms with the same number
of fixed points will have the {\it same}  matrix representation with
respect to an adapted basis, the intersection matrices will not be
the same and, therefore, the corresponding two matrix
representations in the symplectic group will not be conjugate.

%\marginpar{need notion of canonical representative in symplectic
%group}

There is an algorithm to  replace $M_{{\tilde{\mathcal{A}}}}$ by a
the symplectic matrix $M_{h_{CAN}}$ by replacing the submatrix
$N_{(t-2),(p-1)\times(p-1)}$ by a symplectic matrix of the same size
(see \cite{JPSym}). We will call this matrix
$N_{symp{\tilde{\mathcal{A}}}}$  and give an example in  section
\ref{section:ex}.

We note that $I_{{\tilde{\mathcal{A}}}}$ is of the form
$$\left(
\begin{array}{cccc}
 0& I_{pg_0}& 0   & 0\\
-I_{pg_0} & 0 & 0 & 0  \\
0 &0 &B_1 &B_2\\
0& 0& B_3& B_4 \\
\end{array} \right) $$
where the blocks $B_i$ are of the appropriate dimension and we let
$B$ denote the matrix
$$\left(
\begin{array}{cc}
B_1 &B_2\\
 B_3& B_4 \\
\end{array} \right). $$

\section{Schreier-Reidemeister Rewriting} \label{section:SCHR}

If we begin with an arbitrary finitely presented group $G_0$ and a
subgroup $G$, the Schreier-Reidemeister rewriting process tells one
how to obtain a presentation for $G_0$ from the presentation for
$G$. In our case the larger group $G_0$ will correspond to the group
uniformizing $S_0$ and the subgroup $G$ corresponds to the group
uniformizing $S$.

In particular, one chooses a special  set of coset representatives
for $G$ modulo $G_0$, called Schreier representatives, and uses
these to find a set of generators for $G$. These generators are
labeled by the original generators of the group and the coset
representative.

\comment{ these to find in with a finitely presented group $G$ and
find generators and relations for a subgroup $G_0$ given a
homomorphism from the larger group onto a finite group $H$ where the
subgroup is the kernel}
\subsection{The relation between the action of the homeomorphism and
the surface kernel subgroup}
%\subsection{The surface kernel and the rewriting
\label{section:rewrite}

We may assume that $S_0 = U/F_0$ where $F_0$ is the Fuchsian group
with presentation
\begin{equation} \label{equation:presentation}
\langle a_1, ..., a_{g_0}, b_1,..., b_{g_0}, x_1, ...x_t|\;
 x_1 \cdots x_t (\Pi_{i=1}^g[a_i,b_i]) =1;
 x_i^p=1\rangle.
\end{equation}

We summarize the result of \cite{G3} repeating  some facts about the
$n_i$ and $m_j$. We let $\phi: F_0 \rightarrow \mathbb{Z}_p $ be
given by
$$\phi(a_i) = \phi(b_i) =0 \;\;\; \forall i=1,...,g_0 \mbox{ and }
\phi(x_j) = n_j \ne 0 \;\;\;\forall j=1,...,t.$$ The  $n_i$ satisfy
$\Sigma_{i=1}^t n_i \equiv 0 \;(p)$. If $F = Ker\; \phi$, then $S=
U/F$. Moreover, $F_0/F$ acts on $S$ with quotient $S_0$. Conjugation
by $x_1$ acts on $F$ and if $h$ is the induced conformal map on $S$,
$\langle h \rangle$ is isomorphic to the action induced by this
conjugation and the conjugacy class of $h$ in the mapping-class
group is determined by the set of $n_i$. Replacing $h$ by a
conjugate we may assume that $0 < n_i \le  n_j < p$ if $i < j$.

% It is not difficult to show that the order of the $n_i$ does
%not matter.

Since $m_i$ the number of $j$ such that $\phi(x_j) = i$, we have
$\Sigma_{i=1}^{p-1} im_i \equiv 0 \; (p)$ and the conjugacy class of
$h$ is also completely determined by $(m_1,...,m_{p-1})$.

When we need to emphasize the relation of $h$ to $\phi$, we write
$h_{\phi}$ to mean the automorphism determined by conjugation by
$x_1$. The conjugacy class of $h^2$, would then be determined by the
homomorphism $\psi$ with $\psi(x_j) \equiv  2\phi(x_j) \; (p)$ or by
conjugation by $x_1^2$.

%\begin{remark} Put somewhere, This gives an alternate proof that $h_{\phi}$ has
%precisely $2- \tr h_{\phi}$ fixed points \cite{Sah}.
%\end{remark}
%\begin{remark} Need to be careful. We have $h_{\phi}$ has the right rotation angle at its fixed points.
%\end{remark}
We note that in \cite{G3} results are written in greatest possible
generality so that $S_0$ has punctures, some of which are fixed by
the automorphism and others of which are not. Here, we use the
results of \cite{G3} for compact surfaces. $F_0$ is sometimes called
surface kernel and $\phi$ the surface kernel homomorphism
\cite{GPLMS, Harvey, Singerman}.

Any other map from $F_0$ onto $\mathbb{Z}_p$ with the same
$(m_1,....,m_{p-1})$ and with $\phi(x_j) \ne 0 \;\; \forall j$ will
yield an automorphism conjugate to $h$.
\subsection{The rewriting} \label{section:rewrite}

We want to apply the rewriting process to words in the generators of
this  presentation for $F_0$ to obtain a presentation for $F$. We
choose coset representatives for $H = \langle h \rangle$ as
$x_1,x_1^2,x_1^3,...,x_1^p$ and observe that these are a {\sl
Schreier} systems. (see page 93 of \cite{MKS}). That is, that every
initial segment of a representative is again a representative.

The Schreier right coset function assigns to a word $W$ in the
generators of $F_0$, its coset representative $\overline{W}$ and
$\overline{W} = x_1^q$ if $\phi(W)= \phi(x_1^q)$.

If ${\mathfrak{a}}$ is a generator of $F_0$, set
$S_{K,\;{\mathfrak{a}}} = K {\mathfrak{a}}
{\overline{K{\mathfrak{a}}}^{\;-1}}$. The rewriting process $\tau$
assigns to a word that is in the kernel of the map $\phi$, a word
written in the specific generators, $S_{K,\;{\mathfrak{a}}}$ for
$F$.
 Namely,  if ${\mathfrak{a}}_w$, $w = 1, ..., r$ are generators for
 $F_0$ and
$$U= {\mathfrak{a}}_{v_1}^{\epsilon_1}{\mathfrak{a}}_{v_2}^{\epsilon_2} \cdots
{\mathfrak{a}}_{v_r}^{\epsilon_r} \;\;\;\;\;(\epsilon_i = \pm 1),$$
defines an element of $F$, then (corollary 2.7.2 page 90 of
\cite{MKS})
$$\tau(U)= S_{K_1,\;{\mathfrak{a}_{v_1}}}
^{\epsilon_1} S_{K_2,\;{\mathfrak{a}}_{v_2}}^{\epsilon_2} \cdots
S_{K_r,\;{\mathfrak{a}}_{v_r}}^{\epsilon_r}$$ where $K_j$ is the
representative of the initial segment of $U$ preceding
${\mathfrak{a}}_{v_j}$ if $\epsilon_j = 1$ and $K_j$ is the coset
representative of $U$ up to and including
${\mathfrak{a}}_{v_j}^{-1}$ if $\epsilon_j = -1$.

In our case each ${\mathfrak{a}}_v$ stands for some generator of
$F_0$, that is one of the $a_i$ or $b_i$ or $x_j$.

 We apply
theorem 2.8 of \cite{MKS} to see
\begin{theorem} \label{theorem:Fpres} Let  $F_0$ have the presentation
given by equation \ref{equation:presentation}. Then $F$ has presentation
\begin{eqnarray} \label{equation:presentationF}
\lefteqn{\langle S_{K,\;a_i}, \;S_{K,\;b_i},\; i = 1, ..., g_0;\;
S_{K,x_j}, \; j =
1,...t|}&\\
 && \tau(K \cdot x_1 \cdots x_t(\Pi_{i=1}^{g_0}[a_i,b_i]) \cdot  K^{-1})=1,\;  \tau(Kx_j^pK^{-1})=1 \rangle.
\end{eqnarray}
\end{theorem}

\begin{proof}
Let $K$ run over a complete set of coset representative for $\phi:F \rightarrow H$.
Then $F_0$ has generators
\begin{eqnarray*}
  S_{K,\;a_i}, S_{K,\;b_i} & &  i = 1, ..., g_0\\
  S_{K,\;x_j}, && j =1,...t
  \end{eqnarray*}
and relations
\begin{eqnarray} \label{equation:eqpresentationF}
\tau(K(x_1 \cdots x_t\Pi_{i=1}^{g_0}[a_i,b_i])
K^{-1})&=&1\\
 \tau(Kx_j^pK^{-1})&=&1 \end{eqnarray}
\end{proof}

We want to simplify this presentation and eliminate generators and
relations so that there is a single defining relation for the
subgroup. We first assume that $\phi(x_1) = h$. We note that if we
can find a homology basis adapted to $h$, we can easily find a
homology basis adapted to any power of $h$ and, therefore, this
assumption will not be significant.

We will show:

\begin{theorem}\label{theorem:homotopy} Let $F_0$ have the
presentation given by equation \ref{equation:presentation}. Then $F$
has presentation
%\begin{equation}
%\label{equation:presentationF}
$$\langle h^j(A_i), h^j(B_i), i=1,....,g_0, j = 0, ..., p-1:
h^j(X_i), i=3,...,t, j = 0, ..., p-2 | {\hat{\hat{R}}} =1 \rangle.$$
%\end{equation}
The relation $\hat{\hat{R}}$ is the single defining relation for the
group $F$. Each generator and its inverse occur exactly once in
${\hat{\hat{R}}}$ Further, every generator that appears in
${\hat{\hat{R}}}$ is linked to another distinct generator.
\end{theorem} %\marginpar{check terminology}

and
\begin{corollary} The homology basis obtained by abelianizing  the basis in theorem
\ref{theorem:homotopy} gives a homology basis adapted to $h$.
\end{corollary}

We note that an explicit formula for $\hat{\hat{R}}$ is  given in
\cite{JGPSym}.

\begin{proof}

 If we let $\phi(x_1) = h$ and $\phi(K)=\phi(x_1)^r$, then we have
$S_{K,X_j}= x_1^{r\phi(x_1)} \cdot x_j \cdot {\overline{K \cdot
x_j}}^{-1}$. Thus if $X_j = x_j \cdot {\overline{x_j}}^{-1}$, then
$S_{K,x_j} = h^r(X_j)$.

We begin by rewriting the generators and relations using this
notation.

First we find $\tau(\overbrace{x_1x_1\cdots
x_1}^{p-\mbox{factors}})=1$. Setting  $X_1= S_{x_1^p,x_1}$ since
$\overline{1} = x_1^p$, we have
\begin{equation}
\tau(x_1^p) = X_1 \cdot  h(X_1) \cdot h^2(X_1) \cdots
h^{p-2}(X_1)h^{p-1}(X_1)=1.
\end{equation}

Similarly, if $\phi(x_j) = n_j$, and we set $X_j =
S_{{\overline{1}},x_j}$, then if ${\overline{K}} = x_1^{s}$,

then we can write $S_{K,\;x_j} = h^{s\cdot n_j}(X_j)$.

This tells us that
\begin{equation} \label{equation:homotopyXj}
\tau(x_j^p) = X_j \cdot h^{n_j}(X_j) \cdot h^{2n_j}(X_j) \cdots
h^{(p-2)n_j}(X_j) h^{(p-1)n_j}(X_j) =1 . \end{equation}

In particular, we will make special use of this when $j=2$
\begin{equation} \label{equation:homotopyXj2}
\tau(x_2^p) = X_2 \cdot h^{n_2}(X_j) \cdot h^{2n_2}(X_2) \cdots
h^{(p-2)n_2}(X_2) h^{(p-1)n_2}(X_2) =1 . \end{equation}

Note that in deriving all equations we are free to make use of the
fact (see \cite{MKS}) that \begin{equation}
\label{equation:freehomotopy} S_{M,\;x_1} \approx 1 \forall \mbox{
Schreier representatives} M
\end{equation}
 where $\approx$ denotes {\sl freely equal to}.
This eliminates the $p$ generators, $S_{x_1^j, x_1}, j = 1,...,p$. %$\endpf$

Now  equation (\ref{equation:homotopyXj}) is a relation in the
fundamental group. We remind the reader that for a compact Riemann
surface, homology is abelianized homotopy so that when abelianized,
it reduces to
\begin{equation} \label{equation:homologyXj}
h^{p-1}(X_j) {\approx}^h -X_j - h(X_j) - \cdots - h^{p-2}(X_j)
\end{equation}  where
$\approx^h$ denotes {\it is homologous to}.

We also note that the $\tau(Kx_j^pK^{-1})=1$ do not give us any
additional relations for $\bar{K} \ne 1$, but  merely a conjugate
relation already implied by $\tau(x_j^p)=1$.

Next we set $A_i = S_{{\overline{1}},\;a_i}$ and $h^r(A_i) =
S_{x_1^r, \;a_i}$, $B_i = S_{1,\;b_i}$ and $h^r(B_i) = S_{x_1^r,\;
b_i}$, then the relation $R=1$, where $R=x_1 \cdots
x_t(\Pi_{i=1}^{g_0}[a_i,b_i])$, yields $\tau(R) =1$ and we obtain

\begin{equation} \label{equation:tauR}
X_1h^{n_1}(X_2)h^{n_1 \cdot n_2}(X_2) \cdots h^{n_1 \cdot n_2 \cdots
n_{t-2}}(X_{t-1}) \cdot h^{n_1 \cdot n_2 \cdots n_{t-1}}(X_t)
(\Pi_{i=1}^{g_0} [A_i,B_i]) =1
\end{equation}

and using the fact that $X_1 \approx 1$, we have
\begin{equation} \label{equation:tauRX1}
h^{n_1}(X_2)h^{n_1 \cdot n_2}(X_3) \cdots h^{n_1 \cdot n_2 \cdots
n_{t-2}}(X_{t-1}) \cdot h^{n_1 \cdot n_2 \cdots n_{t-1}}(X_t)
(\Pi_{i=1}^{g_0} [A_i,B_i]=1
\end{equation}

Similarly, we obtain the relations $\tau(KRK^{-1}) =1$.

We can solve equation (\ref{equation:tauRX1}) for $(h(X_2))^{-n_1}$
to obtain

\begin{equation} \label{equation:tauRX1X2}
(h(X_2))^{-n_1} = h^{n_1\cdot n_2}(X_3) \cdots h^{n_1 \cdot n_2
\cdots n_{t-2}}(X_{t-1}) \cdot h^{n_1 \cdot n_2 \cdots n_{t-1}}(X_t)
(\Pi_{i=1}^{g_0} [A_i,B_i])
\end{equation}

Now each of the relations $\tau(KRK^{-1})$ allows us to solve for
$h^q(X_2)$ for some $q$ and for each $K$ we obtain a different $q$.
Therefore, we can substitute equation (\ref{equation:tauRX1X2}) and
its images under powers of $h$ into equation
(\ref{equation:homotopyXj2}) (i.e. the relation $\tau(x_2^p) =1$).
We thus eliminate all of the generators of the form $S_{M,x_2}$ for
each coset representative $M$ and all of the relations
$\tau(KRK^{-1})$. We obtain one new relation ${{\hat{R}}}$ from
equation (\ref{equation:homotopyXj2}). This relations involves
$h^w(X_j)$ for every $w = 0,1,..., p-2$ and every $j = 3, ..., t$.
We also note that for each $q$ the sequence
$\Pi_{i=1}^{g_0}[h^q(A_i),h^q(B_i)]$ occurs in ${\hat{R}}$.
 Using equation (\ref{equation:homotopyXj}) to replace the
generator $h^{p-1}(X_j)$ by a word in the $h^{-v}(X_j)$ with $v = 0,
..., p-1$ we eliminate those relations and we obtain a single
defining relation $ {\hat{\hat{R}}}$ involving each of the following
generators below  and their inverses exactly once. The generators
are
$$h^d(X_j) \;\; j= 3,...,t \mbox{ and } d = 0, 1, ..., p-2$$ and $$
h^d(A_i), h^d(B_i) \;\; i= 1,..., g_0 \mbox{ and } d = 0,1, ...,
p-1$$

We obtain  $2g_0  p  + (t-2)(p-1)$ generators and a single defining
relation. It is fairly straight forward to check that the relation
has the last two properties of the theorem.

\end{proof}

We recapitulate. The idea of this proof is that for each $K$, we can
solve $\tau(KRK^{-1})= 1$ for an appropriate image $S_{K',x_2}$
where $K'$ depends upon $K$. The appropriate image of $S_{K',x_2}$
is placed on the left of the equation, and we then substitute the
right hand side of the solution into the relation $\tau(x_2^p) =1$.
This yields one relation ${\hat{{\mathcal{R}}}}$ which involves each
$A_i,B_i$ generator and their inverses and all of their images under
powers of $h$ and each $S_{K,x_j}$ and the images under powers of
$h$ but no inverses. We still have the finite order relations
$\tau(x_j^p)$, $j = 3, ...., t$. The $\tau(Kx_j^pK^{-1})=1$ are
merely permutations of the relation $\tau(x_j^p) =1$ so we can
eliminate all but one of these. For each $j = 3, ..., t, \;\;
\tau(x_j^p)$ can be solved for $h^{p-1}(S_{{\overline{1}},x_j})$. It
will be a word in the inverses of all of the other $S_{K,x_j}$. We
substitute these into ${\hat{{\mathcal{R}}}}$ and obtain a single
defining relation ${\hat{\hat{R}}}$ in which every generator and its
inverse occurs exactly once.

\vskip .1in
 We introduce further
terminology.

\begin{definition} A relation is {\sl evenly worded} if for each generator  $A$
that occurs in the relation $A^{-1}$ also occurs. The generators $A$
and $B$ occurring in a relation are {\sl  linked} if  the relation
is of the form $W_0AW_1BW_2A^{-1}W_3B^{-1}W_4$ where the $W_i,
i=0,...,4$ are words in the generators not involving $A^{\pm1} $ or
$B^{\pm 1}$. The relation is {\sl fully linked} if each generator
$A$ occurring in the relation is linked to a unique distinct
generator.
\end{definition}

We can,  therefore, say that ${\hat{\hat{R}}}$ is evenly worded and
fully linked.
\section{The Case $t=0$} \label{section:t=0}

If the number of fixed points is zero, we can still find an adapted
basis. The calculations are slightly different. We have $2g =
2p(g_0-1) + 2$.  The presentation given by
(\ref{equation:presentation}) for the group $F_0$ becomes

\begin{equation} \label{equation:0presentation}
\langle a_1, ..., a_{g_0}, b_1,..., b_{g_0} |\;
  (\Pi_{i=1}^{g_0}[a_i,b_i]) =1;
 \rangle.
\end{equation}

Again, replacing $h$ by a conjugate if necessary, the map $\phi: F_0
\rightarrow \mathbb{Z}_p$ can be taken to be
$$\phi(a_i) = \phi(b_i) =0 \;\;\; \forall i=2,...,g_0 \mbox{ and }
\phi(a_1) = 1 \mbox{ and  } \phi(b_1) =0 $$

Using the rewriting with coset representatives
$1,a_1,...,a_1^{p-1}$, note that $S_{a_1^{k},a_1} \approx 1$ for $k
= 0, ... ,p-2$. Let $A = S_{a_1^{p-1},a_1}$. Then $h$ acts on $Ker
\; \phi$ via conjugation by $a_1$. We have $h(A) = A$

We let $\{h^k(A_j),h^k(B_j), j = 2,...,g_0, k=0,...,p-1\}$ be as in
the proof of Theorem \ref{theorem:homotopy} and let $B= S_{1,b_1}$.

%Notice that  $h^k(A_j) = \tau(a_1^k)$. \marginpar{why?}

 We let $P = \Pi_{i=2}^{g_0}
[A_i,B_i]$.

Then we can compute that
\begin{align} \label{equation:new}
\tau(R) = 1 \implies h(B)B^{-1}T =1\\
\nonumber \tau(a_1^kRa_1^{-k}) = 1 \implies
h^k(B)(h^{k-1}(B))^{-1}h^k(P)
= 1 \forall k=1,...p-2\\
\nonumber  \mbox{and}\; \tau(a_1^{p-1}Ra_1^{-(p-1)})=1 \implies
ABA^{-1}(h^{p-1}(B))^{-1}h^{p-1}(P) = 1
\end{align}

We eliminate generators using (\ref{equation:new}) and let $\alpha =
A$ and $\beta = h^{p-1}(B)$ to obtain generators
$$\{\alpha, \beta \} \cup \{h^k(A_j),h^k(B_j), j = 2,...,g_0,
k=0,...,p-1\}$$ and the single defining relation $\beta \alpha
\beta^{-1}= h^{p-1}(P) \alpha \Pi_{i=k}^{p-2}h^k(P).$ Further we
calculate that $h(\alpha) \approx^h \alpha$ and $h(\beta) \approx^h
\beta$. Note that $P$ is a product of commutators. Thus the matrix
representation for $h$ on $S= U/F$ where $F = \mbox{Ker} \phi$ is
given by $2(g_0-1)$ permutation matrices $M_{p \times p}$ and one
two by two identity matrix.

The basis is a canonical homology basis. We have $\alpha \times
\beta =1$. Further $\alpha$ and $\beta$ are disjoint from any other
lifted curves in the basis. Each lift of a generator other than
$a_1$ or $b_1$ is a simple closed curve and is disjoint from all
other curves except for the corresponding $h^k(B_j)$. That is,
$h^k(A_i) \times h^r(B_j) = \delta_{ij} \cdot \delta_{kr}$ for all
integers $i,j \in \{2,...,p-1\}$ and $k,r \in \{0,...,p-1\}$. We
replace $\alpha$ or $\beta$ by an appropriate conjugate if
necessary. We note that the curves $\alpha$ and $\beta$ are by
default of type (1) in definition \ref{definition:adapt}

We have obtained the  following version of Theorem \ref{theorem:int}
when $t=0$
\begin{theorem} {\label{theorem:0int}}
If $t=0$, then the surface $S$ has a canonical homology basis
consisting of:
\begin{enumerate}
\item  $h^j(A_w), h^j(B_w)$  where $2 \le w \le g_0,  0 \le j \le p-1$.

\item $\alpha, \beta$
  where $h^k(\alpha) \approx^h \alpha, h^k(\beta) \approx^h \beta, 0 \le k \le p$

\item
 The intersection numbers for the elements of the
adapted basis are given by
\begin{enumerate}
\item $h^j(A_w) \times h^j(B_w) = 1$

\item $\alpha \times \beta =1$
\item  All other intersection numbers are $0$ except for those that
following from the above by applying the identities below to
arbitrary homology classes $C$ and $D$.
 \subitem $ C \times D = - D
\times C$ \subitem $ h^j(C) \times h^k(D)  = h^0(C) \times
h^{k-j}(D)$, ($k-j$ reduced modulo $p$.)
\end{enumerate}
\end{enumerate}
\end{theorem}

\section{Example, $p=3$, $t=5$, $(1,1,2,1,1)$} \label{section:ex}
In this section we work out the specific example with $p=3$ and
$t=5$. Assume   $\phi(x_1) = 1 ,\phi(x_2) = 1,\phi(x_3) =
2,\phi(x_4) = 1,\phi(x_5) = 1$ so that $(n_1, ..., n_5) =
(1,1,2,1,1)$ and
 $(m_1,m_2)= (4,1)$.

 First replacing  $h$ by a conjugate,  we may assume that
 $\phi(x_1) = 1 ,\phi(x_2) = 1,\phi(x_3) =
1,\phi(x_4) = 1,\phi(x_5) = 2$. We choose as coset representatives
$x_1, x_1^2$ and $x_1^3$.

For any $g_0$, we have generators $$h^q(A_i),h^q(B_i), q= 0, \ldots,
p-1=2, i=1 \ldots, g_0.$$ and $$S_{{\overline{x_1^r}},x_j}, \;\;\; r
= 1,2,3, j = 2,3,4.  $$ We also have by equation
(\ref{equation:freehomotopy})
$$S_{{\overline{x_1^r}},x_1} \approx 1, r = 1,2,3.$$
and, therefore,  these generators and the relation $\tau(x_1^3)$
drops out of the set of generators and relations to be considered.

We have
\begin{equation}
\tau(x_j^3) = S_{{\overline{1}}, x_j} \cdot S_{{\overline{x_j}},x_j}
\cdot S_{{\overline{x_j^2}},x_j}
\end{equation}

Set   $S_{{\overline{1}}, x_1}= Y_1$. Then $h(Y_1) = S_{x_1,x_1}$
and
 $h^2(X_j) = S_{x_1^2,x_j}$;

Set $S_{{\overline{x_1}}, x_2}= Y_2$ Then $h(Y_2) = S_{x_1^2,x_2}$
and
 $h^2(Y_2) = S_{x_1^3,x_2}$;

 Similarly, set
$S_{{\overline{x_1x_2}}, x_3}= Y_3$.  Then $h(Y_3) = S_{x_1^3,x_3}$
and
 $h^2(Y_3) = S_{x_1,x_3}$;
If $S_{{\overline{x_1x_2x_3}}, x_4}= Y_4$ Then $h(Y_4) =
S_{x_1,x_4}$ and
 $h^2(Y_4) = S_{x_1^2,x_4}$;and
finally if
 $S_{{\overline{x_1x_2x_3x_4}},
x_5}= Y_5$, then $h(Y_5) = S_{x_1^2,x_5}$ and
 $h^2(Y_5) = S_{x_1^3,x_j}$;

Use this notation and use $\tau(x_j^3)=1$ to see that

\begin{align} \label{equation:squares}
  %to remove numbering (before each equation)
h^2(Y_2)\cdot Y_2 \cdot h(Y_2)= 1 \\
\nonumber h(Y_3)\cdot h^2(Y_3)\cdot Y_3=
1 \\
\nonumber Y_4\cdot h(Y_4)\cdot h^2(Y_4) = 1\\ \nonumber
h^2(Y_5)\cdot h(Y_5)\cdot Y_5 = 1 \end{align}

We compute \begin{equation} \label{equation:tauR} \tau(R) =
  S_{{\overline{1}},x_1}\cdot
 S_{{\overline{x_1}},x_2} \cdot
S_{{\overline{x_1x_2}},x_3} \cdot S_{{\overline{x_1x_2x_3},x_4}}
\cdot S_{{\overline{x_1x_2x_3x_4}},x_5} \cdot (\Pi_{i=1}^{g_0}
[A_i,B_i])=1.
\end{equation}
Using $S_{{\overline{x_1^r}},x_1} \approx 1, r = 1,2,3$ and solving
for $(Y_2)^{-1}$ in \ref{equation:tauR}, we have
\begin{equation}
(S_{x_1,x_2})^{-1} = Y_2^{-1} = Y_3 \cdot Y_4 \cdot Y_5  \cdot
(\Pi_{i=1}^{g_0} [A_i,B_i])=1.
\end{equation}
and
\begin{equation}
 (h(Y_2))^{-1} = h(Y_3) \cdot h(Y_4) \cdot h(Y_5)  \cdot
(\Pi_{i=1}^{g_0} [h(A_i),h(B_i)])=1.
\end{equation}
\begin{equation}
  (h^2(Y_2))^{-1} = h^2(Y_3) \cdot h^2(Y_4) \cdot h^2(Y_5)  \cdot
(\Pi_{i=1}^{g_0} [h^2(A_i),h^2(B_i)])=1.
\end{equation}
Using  $h^2(Y_2)Y_2h(Y_2) = 1$ and letting  $P=(\Pi_{i=1}^{g_0}
[A_i,B_i])$, we have
\begin{equation}
h(Y_3) \cdot h(Y_4) \cdot h(Y_5)  \cdot h(P)  \cdot Y_3 \cdot Y_4
\cdot Y_5  \cdot P) \cdot h^2(Y_3) \cdot h^2(Y_4) \cdot h^2(Y_5)
\cdot h^2(P) =1
\end{equation}
We use equation  (\ref{equation:squares})  to replace the $h^2(Y_j)$
and
 obtain the relation
% equation \ref{align:squares}?
%\begin{equation} \label{equation:hathatR}
%\end{equation}
\begin{align} \label{equation:hathatR}
  %to remove numbering (before each equation)
  h(Y_3) \cdot h(Y_4)
 \cdot h(Y_5)  \cdot h(P) &  \cdot Y_3 \cdot Y_4
\cdot Y_5  \cdot P) \cdot (h(Y_3))^{-1} \cdot (Y_3)^{-1}\\
 \nonumber &\cdot (h(Y_4))^{-1} \cdot (Y_4)^{-1}   \cdot (Y_5)^{-1}
\cdot (h(Y_5))^{-1}\cdot h^2(P)=1.
\end{align}

Equation (\ref{equation:hathatR}) is the relation
${\hat{\hat{R}}}=1$.

We use the notation of section \ref{section:matrixforms} in
particular the defintion of the matrix $B$ given at the end of that
section.  We can compute from the formulas for intersection numbers
in Theorem \ref{theorem:GPint} that the relevant part of the
intersection matrix, $B$, is  the $6 \times 6$ submatrix that gives
the intersection matrix for the curves in the basis given in the
order
$$X_{1,3}, h(X_{1,3}),X_{1,4},h(X_{1,4}),X_{2,1}, h(X_{2,1})$$ is
$$B= I_{\hat{\hat{R}}}= \left(
  \begin{array}{rrrrrr}
    0 & 1 & 1 & 0 & 1 &-1 \\
    -1 & 0 &-1 & 1 & 0 & 1 \\
    -1 & 1 & 0 & 1 & 1 & -1 \\
    0 & -1 & -1 & 0 & 0 & 1 \\
    -1 & 0 & -1 & 0& 0 & 0 \\
    1 & -1 & 1 & -1 & 0 & 0 \\
  \end{array}
\right).$$

That is, the matrix $I_{{\tilde{\mathcal{A}}}}$ breaks up into
blocks $\left(
  \begin{array}{ccc}
    0 & I_{pg_0} & 0 \\
    -I_{pg_0} & 0 & 0 \\
    0 & 0 & I_{\hat{\hat{R}}} \\
  \end{array}
\right)$

 We
now rearrange the relation. To simply the notation we let $a= Y_3$,
$b=Y_4$ and $c = Y_5$. So that the relation becomes
\begin{equation}
h(a)\cdot h(b) \cdot h(c) \cdot h(P) \cdot a \cdot b \cdot c \cdot P
\cdot (h(a))^{-1}a^{-1}(h(b))^{-1} \cdot b^{-1} \cdot c^{-1} \cdot
(h(c))^{-1}\cdot h^2(P) = 1. \end{equation}

We can also make the simplifying assumption, replacing the elements
that occur in $P$, $h(P)$ and $h^2(P)$ by conjugates,  that we are
merely working with the symbol

\begin{equation}
h(a)\cdot h(b) \cdot h(c)  \cdot a \cdot b \cdot c  \cdot
(h(a))^{-1}a^{-1}(h(b))^{-1} \cdot b^{-1} \cdot c^{-1} \cdot
(h(c))^{-1} = 1. \end{equation}

We replace generators and relations using the algorithm of
\cite{JGSym}  as follows:

Let $M=h(a)\cdot W_1\cdot h(b) W_2\cdot (h(a))^{-1}$ where $W_1=
\emptyset$, $W_2 = h(c) \cdot a \cdot b \cdot c$. Set $W_3 = a^{-1}$
and $W_4 = b^{-1} \cdot c^{-1} \cdot (h(c))^{-1}$.

Let $N= W_3W_2(h(a))^{-1}$. Then

\begin{align} \label{equation:rel}
  h(a)\cdot h(b) \cdot h(c) \cdot a \cdot b
\cdot c  \cdot (h(a))^{-1}a^{-1}(h(b))^{-1} \cdot b^{-1} \cdot
c^{-1}\cdot (h(c))^{-1} \\
\nonumber  = [M,N]W_3W_2W_1W_4 \;\;\;\;\;\;\;\;\;\; \;\;\;\;\;\;\;\;\;\;\;\;\;\;\;\\
\nonumber
 = [M,N] \cdot
a^{-1} \cdot h(c) a b c b^{-1} \cdot c^{-1} \cdot (h(c))^{-1}\\
\nonumber =1.
\end{align}

At this point one can proceed by inspection and let  $[b,c]^*$
denote the conjugate of $[b,c]$ by $a^{-1} \cdot h(c) \cdot a$ to
obtain

$$[M,N]\cdot [b,c]^* \cdot [a^{-1}, h(c)] =1. $$

However, to follow the algorithm carefully, we would set $\tilde{M}
= a^{-1}\cdot h(c) \cdot a$ and $\tilde{N}= b \cdot c \cdot b^{-1}
\cdot c^{-1} \cdot a$.

Then equation (\ref{equation:rel})  becomes $$[M,N]\cdot [\tilde{M},
\tilde{N}] \cdot [b,c] =1.$$

Thus the canonical homology basis is given by $$\{ h^j(A_i),
h^j(B_i) \}, i= 1,...,g_0, j = 0, ..., p-1 \cup \{M, N, \tilde{M},
\tilde{N}, b, c \}.$$

The reordered basis
 $\{M,
\tilde{M}, b, N, \tilde{N}, c \}$ has intersection matrix
$$\left(
                                               \begin{array}{rrrrrr}
                                                 0  & 0 & 0 & 1 & 0 & 0 \\
                                                 0 & 0 & 0 & 0 & 1 & 0 \\
                                                 0 & 0& 0& 0& 0& 1\\
                                                 -1 & 0& 0& 0& 0& 0\\
                                                 0 & -1 & 0& 0& 0& 0\\
                                                 0 & 0& -1& 0& 0& 0\\
                                                                                              \end{array}
                                             \right). $$

 We
can compute the action of $h$ on these last six elements of the
homology basis. First we note that $h(b) \approx^h M -  c - b -a -
h(c)$ and $h(a) \approx^h -N  + h(c)+ b + c$. Therefore,

$c \mapsto h(c) \approx^h \tilde{M}$

$h(c) \mapsto -c - h(c) \approx^h -c - \tilde{M}$

$a \mapsto h(a) \approx^h -N + h(c) + b +c  \approx^h -N + \tilde{M}
+ b + c$

$b \mapsto h(b) \approx^h M -c -b -a -h(c)\approx^h M -c - b
-{\tilde{N}} - \tilde{M}$

$M \mapsto h(M) \approx^h  -\tilde{M} - N$, and

$N \mapsto h(N) \approx^h -c + M -N $.

Thus the matrix of the action of $h$ with respect to the ordered
basis $M,\tilde{M},b,N, \tilde{N}, c$ is the submatrix we have been
seeking. Namely,

$$N_{symp{\tilde{\mathcal{A}}}} = \left(
  \begin{array}{rrrrrr}
    0&1&0&-1&0&0 \\
    0 & -1 & 0&1&0&-1 \\
    1&-1&-1&-0&-1&-1 \\
    1&0&0&-1&0& -1\\
    0&1&1&-1&0&-1 \\
    0&1&0&0&0&0 \\
 \end{array}
 \right)$$

One can verify that this $6 \times 6$ matrix really is a submatrix
of a symplectic matrix, as it should be.

\section{Remarks} \label{section:additional}

Recall (section \ref{section:tor}) that there is a surjective map
$\pi: MCG(S_g) \rightarrow SP(2g,\mathbb{Z})$.
 It is well known that the restriction of $\pi$ to elements (mapping-classes) of finite order
 is an isomorphism. It  is shown in
\cite{G4} that Theorem \ref{theorem:ad} implies a stronger result
than this which we note for completeness.

\begin{corollary} \cite{G4}
 If $h$ is a conformal automorphism of $S$ of genus $g \ge 2$ and if there
are two pairs of curves $C_1$, $D_1$ and $C_2$, $D_2$ with $C_i
\times C_j = D_i\times D_j= 0$ for $i=1,2$ and $j= 1,2$ and $C_i
\times Dj = \delta_{ij}$ where $\delta_{ij}$ is the Kronecker delta,
and $h(C_i) \approx^h C_i$ and $h(D_i) \approx^h D_i$ for $i=1,2$,
then $h$ is the identity.
\end{corollary}

\begin{proof} If $h$ is of prime order, simply write each of the four curves as a sum of the
curves in the adapted homology basis, apply $h$ and equate
coefficients. If $h$ is not of prime order, apply this to every
power that is of prime order to see that each must be the identity.
\end{proof}

We obtain immediately,

\begin{corollary} If two mapping-classes of finite order have the
same action on homology, then they are equal. Equivalently, the
restriction of $\pi$ to elements of finite order is an isomorphism
\end{corollary}

 The idea of an adapted homology basis predates Thurston's notation
of a reducible mapping-class. However, it is clear that the notions
are related. A homeomorphism  $h$ is {\sl reducible} if  $h$ fixes a
{\sl a partition} on the surface, that is,  a set of disjoint simple
closed curves on the surface. A mapping-class is reducible if it
contains a reducible representative.

\begin{corollary} If $h$ represents a mapping-class of prime order,
$p\ge 2$,   and $g_0 \ne 0$, then $h$ is a reducible mapping-class.
\end{corollary}

\begin{proof}  Replace $h$ by a conformal representative if necessary.
An element of type $(1)$  (definition \ref{definition:adapt}) in an
adapted basis for $h$ taken along with its images gives a set of
closed curves on the surface, fixed by the homeomorphism of the
surface. Once we have shown that these are simple closed curves (as
we do in theorem \ref{theorem:GPint}), is clear that the element $h$
is reducible, that is, $h$   fixes a {\sl a partition},  a set of
disjoint simple closed curves on the surface.
\end{proof}

\section{Acknowledgements} The author thanks Yair Minsky and the Yale Mathematics
Department for their hospitality and support while some of this work
was carried out.

\vskip .1in  We give a long, but by no means exhaustive
bibliography.

\vskip .1in

Mathematics Department

Rutgers University

Newark, NJ 07102

e-mail: gilman@andromeda.rutgers.edu; jane.gilman@yale.edu
\end{document}